\documentclass[11pt]{article}


\input epsf

\usepackage{amssymb}
\usepackage{amsmath}
\usepackage{epsfig}

\def\tfrac#1#2{{\frac{\lower.6ex
\hbox{$\scriptstyle#1$}}
{\raise.7ex
\hbox{$\scriptstyle#2$}}}}

\def\RR{{\mathbb R}}
\def\CC{{\mathbb C}}

\def\dsp#1{\displaystyle#1}

\def\F11#1#2#3{
{}_1F_1\left(
\begin{array}{c}
\begin{array}{c}\hskip-10pt#1\end{array}\\
\begin{array}{c}\hskip-10pt #2\end{array}
\end{array}
\hskip-8pt;\,#3
\right)}

\begin{document}

\title{New Series Expansions of the Gauss Hypergeometric Function}

\author{Jos\'e Luis L\'opez \\
    Departamento de Ingenieria Mat\'ematica e Inform\'atica,\\
    Universidad P\'ublica de Navarra, 31006-Pamplona, Spain\\
       \and
    Nico M. Temme\\
     IAA, 1391 VD 18, Abcoude,
The Netherlands\footnote{Former address: CWI, Science Park 123, 1098 XG, Amsterdam, The Netherlands} \\
  \\   { \small e-mail: {\tt
    jl.lopez@unavarra.es, 
    nicot@cwi.nl}}
    }

\date{\ }
\maketitle

\begin{abstract}
\noindent
The Gauss hypergeometric function ${}_2F_1(a,b,c;z)$ can be computed by using the power series in powers of $z, z/(z-1), 1-z, 1/z, 1/(1-z),(z-1)/z$. With these expansions ${}_2F_1(a,b,c;z)$ is not completely computable for all complex values of $z$. As pointed out in  Gil, {\it et al.} [2007, \S2.3], the points  $z=e^{\pm i\pi/3}$ are always excluded from the domains of convergence of these expansions.  B\"uhring [1987] has given a power series expansion that allows computation at and near these points. But, when $b-a$ is an integer, the coefficients of that expansion become indeterminate and its computation requires a nontrivial limiting process. Moreover, the convergence becomes slower and slower in that case. In this paper we obtain new expansions of the Gauss hypergeometric function in terms of rational functions of $z$ for which the points $z=e^{\pm i\pi/3}$ are well inside their domains of convergence . In addition,  these expansion are well defined when $b-a$ is an integer and no limits are needed in that case. Numerical computations show that these expansions converge faster than B\"uhring's expansion for $z$ in the neighborhood of the points $e^{\pm i\pi/3}$, especially when $b-a$ is close to an integer number.

\end{abstract}

\vskip 0.8cm \noindent
{\small
2000 Mathematics Subject Classification:
33C05; 41A58; 41A20, 65D20.
\par\noindent
Keywords \& Phrases:
Gauss hypergeometric function. Approximation by rational functions. Two and three-point Taylor expansions.
}
\section{Introduction}\label{sec:intro}

The power series of the Gauss hypergeometric function ${}_2F_1(a,b,c;z)$,
\begin{equation}\label{definition}
{}_2F_1(a,b,c;z)=\sum_{n=0}^\infty{(a)_n(b)_n\over (c)_n n!}z^n,
\end{equation}
converges inside the unit disk. For numerical computations we can use the right hand side of \eqref{definition} to compute ${}_2F_1(a,b,c;z)$ only in the disk $\vert z\vert\le\rho<1$, with $\rho$ depending on numerical requirements, such as precision and efficiency. From \cite[\S\S 2.3.1 and 2.3.2]{gil} or \cite[eq. 15.2.1 and \S\S 15.8(i) and 15.8(ii)]{nistgauss} we see that the Gauss hypergeometric function ${}_2F_1(a,b,c;z)$ may be written in terms of one or two other ${}_2F_1$ functions with any of the following arguments

\begin{equation}\label{rational}
{1\over z}, \hskip 1cm 1-z,\hskip 1cm {1\over 1-z},\hskip 1cm {z\over 1-z},\hskip 1cm {z-1\over z}.
\end{equation}
As explained in \cite[\S2.3.2]{gil}, when these formulas are combined with the series expansion \eqref{definition}, we obtain a set of series expansions of ${}_2F_1(a,b,c;z)$ in powers of some of the rational functions given in \eqref{rational}. The domains of convergence of the whole set of the expansions obtained in this way are the regions
\begin{equation}\label{regions}
\begin{array}{lll}
\dsp{
\left\vert z\right\vert\le\rho<1,} \quad &
\dsp{\left\vert{1\over z}\right\vert\le\rho<1, }\quad &
\dsp{\vert1-z\vert\le\rho<1,}\\[8pt] 
\dsp{\left\vert{1\over 1-z}\right\vert\le\rho<1,}\quad &
\dsp{\left\vert{z\over 1-z}\right\vert\le\rho<1,}\quad &
\dsp{\left\vert{z-1\over z}\right\vert\le\rho<1.}
\end{array}
\end{equation}

\bigskip
%
\centerline{\includegraphics[width=2.05in]{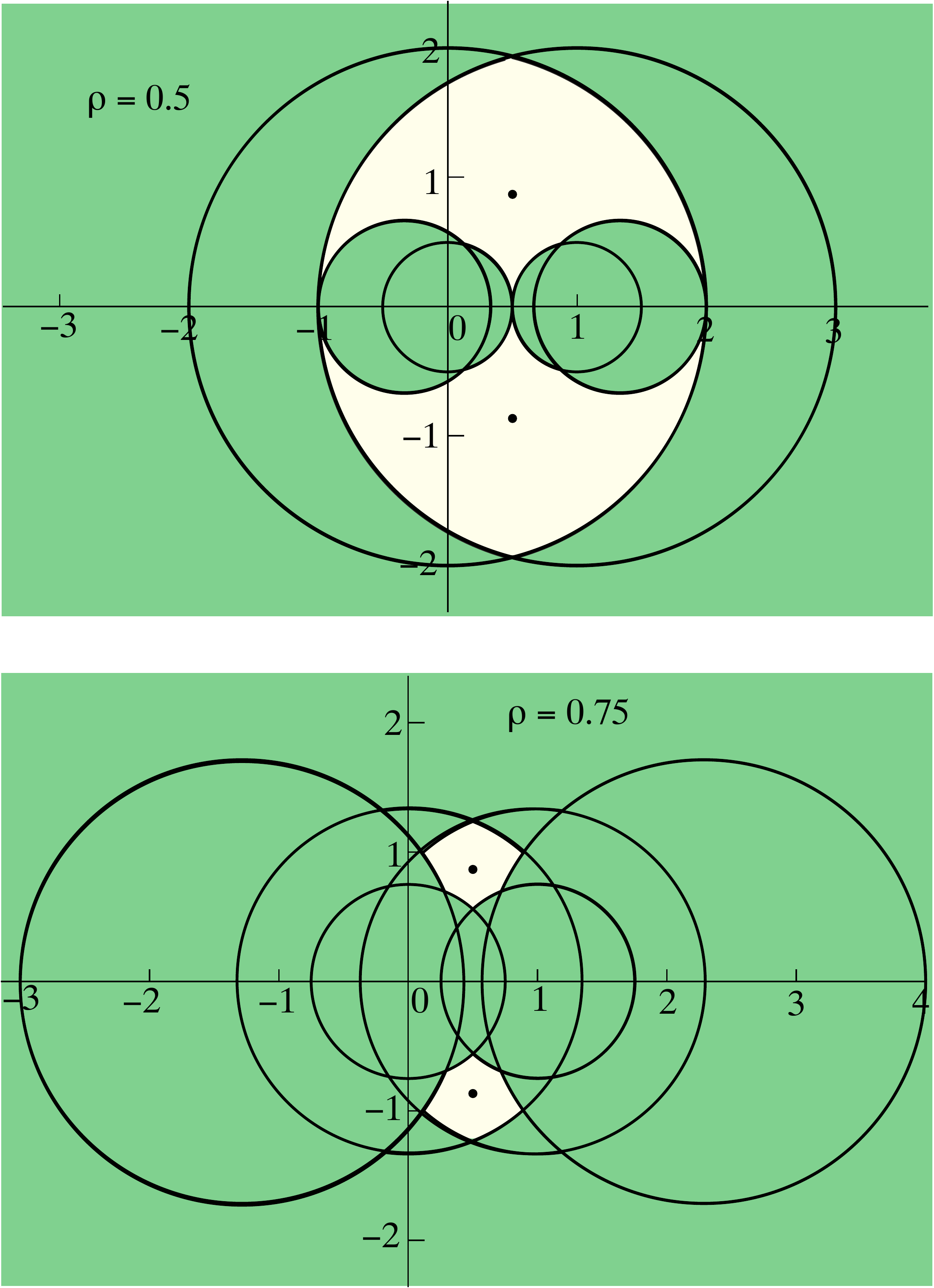}}
%
\noindent
{\bf Figure 1.} Regions given by formulas \eqref{regions} for two particular values of $\rho$. The points $e^{\pm i\pi/3}$ are indicated by black dots.
\medskip

These regions (interior or exterior of certain circles) do not cover the entire $z-$plane, the points $z=e^{\pm i\pi/3}$, that are the intersection points of the circles $\vert z\vert=1$ and $\vert z-1\vert=1$, are excluded  for any value of $\rho<1$ (see Fig. 1).

When $\rho\to 1$, the set of points of the $z-$plane excluded from the union of these regions shrinks to the exceptional points $z=e^{\pm i\pi/3}$, in addition, the convergence of those expansions becomes slower and slower when $z\to e^{\pm i\pi/3}$. To compute the Gauss hypergeometric function in a neighborhood of these points, other methods are indicated in \cite{gil}, the most useful one being B\"uhring's analytic continuation formula \cite{buring}. B\"uhring's expansion reads as follows \cite[\S2.3.2]{gil}. If $b-a$ is not an integer, we have for $\vert{\rm ph}(z_0-z)\vert<\pi$ the continuation formula
\begin{equation}\label{bur}
\begin{array}{ll}
{}_2F_1(a,b,c;z)= & {\Gamma(c)\Gamma(b-a)\over\Gamma(b)\Gamma(c-a)}(z_0-z)^{-a}\sum_{n=0}^\infty d_n(a,z_0)(z-z_0)^{-n}+ \\ & \hskip 2cm {\Gamma(c)\Gamma(a-b)\over\Gamma(a)\Gamma(c-b)}(z_0-z)^{-b}\sum_{n=0}^\infty d_n(b,z_0)(z-z_0)^{-n}, 
\end{array}
\end{equation}
where both series converge outside the circle $\vert z-z_0\vert=\max\lbrace\vert z_0\vert,\vert z_0-1\vert\rbrace$ and the coefficients are given by the three-term recurrence relation
$$
\begin{array}{ll}
d_n(s,z_0)=&{n+s-1\over n(n+2s-a-b)}\lbrace z_0(1-z_0)(n+s-2)d_{n-2}(s,z_0)+[(n+s)(1-2z_0)\\
&+(a+b+1)z_0-c]d_{n-1}(s,z_0)\rbrace,
\end{array}
$$
with $n=1,2,3,...$ and starting values 
$$
d_{-1}(s,z_0)=0, \hskip 2cm d_0(s,z_0)=1.
$$

For the case when $b-a$ is an integer, the coefficients of expansion \eqref{bur} become indeterminate and a limiting process is needed (see \cite{buring} for further details). When we take $z_0=1/2$, the series in \eqref{bur} converges outside the circle $\vert z-1/2\vert=1/2$, and both points $z=e^{\pm i\pi/3}$ are inside the domain of convergence. But, when $b-a$ approaches an integer value, the convergence of the expansion becomes slower and slower.

In this paper we investigate new convergent expansions of the Gauss hypergeometric function ${}_2F_1(a,b,c;z)$ that include the points $z=e^{\pm i\pi/3}$ inside their domain of convergence and which do not require any further computation when $b-a$ is an integer. The starting point is the integral representation \cite[eq. 15.6.1]{nistgauss}
\begin{equation}\label{integraluno}
{}_2F_1(a,b,c;z)={\Gamma(c)\over\Gamma(b)\Gamma(c-b)}\int_0^1t^{b-1}(1-t)^{c-b-1}(1-zt)^{-a}dt,
\end{equation}
valid for $\Re c>\Re b>0$ and $\vert{\rm ph}(1-z)\vert<\pi$.

When we replace $f(t):=(1-zt)^{-a}$ in this integral by the standard Taylor series expansion of $f(t)$ at $t=0$ and interchange summation and integration, we obtain the power series expansion \eqref{definition}. The Taylor series expansion of $f(t)$ at $t=0$ converges uniformly in $z$ for any $t\in(0,1)$ (for any $t$ in the integration domain of \eqref{integraluno}) if $\vert z\vert<1$. Then, the expansion \eqref{definition} is convergent in the disk $\vert z\vert<1$. 

For purposes that will become clear later, it is more convenient to consider the above argument about the region of convergence of the right hand side of \eqref{definition} from a different point of view, which is the following. The domain of convergence of \eqref{definition} (the disk $\vert z\vert<1$) is determined by the two following requirements: (i) The interval of integration $(0,1)$ in \eqref{integraluno} must be completely contained in the domain of convergence of the series expansion of $f(t)$, a disk $D_r$ of center $t=0$ and radius $r\ge 1$, $D_r=\lbrace t\in\CC$, $\vert t\vert<r\rbrace$. (ii) The branch point $t=1/z$ of $f(t)$ must be located outside that domain $D_r$, which means that $z$ must be located in a region $S_r=$ {\it the inverse to the exterior of} $D_r$: $S_r=(D_r^{\rm EXT})^{-1}=\lbrace z\in\CC$, $\vert z\vert<r^{-1}\rbrace$. Therefore, the smaller $D_r$ is (the smaller $r$), the bigger the domain $S_r$ of validity of \eqref{definition} is. But $D_r$ must satisfy $(i)$ and then the largest possible $r$ is $r=1$ and  $S_r=\lbrace z\in\CC$, $\vert z\vert<1\rbrace$ (see Fig. 2).

In this paper we explore the following idea. Instead of the Taylor series expansion of $f(t)$ at $t=0$, consider new different convergent expansions of $f(t)$ in a certain domain $D$ satisfying the two above mentioned requirements: 

\noindent
(i) $(0,1)\subset D$ (The  interval of  integration $(0,1)$ must be completely contained in $D$);

\noindent
(ii) $z\in S:=(D^{\rm EXT})^{-1}$ ($z$ must be located in a region $S=$ {\it the inverse to the exterior of} $D$).

Then, replacing $f(t)$ in \eqref{integraluno} by this new expansion and interchanging summation and integration, we will obtain an expansion of ${}_2F_1(a,b,c;z)$ convergent for $z\in S$. The larger $S$ is, the better, and one expects that, the smaller $D$ is (containing the interval $(0,1)$ in its interior), the bigger $S$ will be. The first possibility that we explore in Section 2 is an expansion of $f(t)$ at $t=1/2$, halfway the  interval of integration $(0,1)$. In Section 3 we generalize this idea expanding $f(t)$ at a generic point $t=w$. In Section 4 we explore a two-point Taylor expansion of $f(t)$ at $t=0$ and $t=1$.
In Section 5 we explore a three-point Taylor expansion of $f(t)$ at $t=0$, $t=1/2$ and $t=1$.
Some final remarks and comments are given in Section 6.

Before we conclude this section, we want to mention other methods for computing the Gauss hypergeometric function. Continued fractions \cite[Chap.~15, \S3]{cuyt},  and Pad\'e approximation \cite{socotres}, \cite[Chap.~9, \S2.4]{gil} , \cite[Chap.~15, Sec. 4]{cuyt},give uniformly convergent expansions on compact subsets of $\CC\setminus[1,\infty)$. In both methods, the approximation is only known explicitly for the exceptional case $a=1$ and $c>b>0$. A different approach using optimal conformal mappings and re-expansions is considered in \cite{socotres} to approximate the Gauss hypergeometric function for $\vert z\vert\ge 1$. A regularization and re-expansion method is used in \cite{socouno} and \cite{socodos} for computing $F(a, b, c; z)$ in the neighborhoods of singular points $z = 1$ and $z = \infty$.

\section{An expansion for $\Re z<1$}

Consider the Taylor expansion of the function $f(t)=(1-zt)^{-a}$ at $t=1/2$:
\begin{equation}\label{tayloruno}
f(t)=\sum_{n=0}^\infty{z^n(a)_n\over n!}\left(1-{z\over 2}\right)^{-a-n}\left(t-{1\over 2}\right)^n.
\end{equation}
This expansion satisfies condition (i) for $D=\lbrace t\in\CC$, $\vert t-1/2\vert<1/2\rbrace$ and $(0,1)\subset D$ (the disk $D$ is the minimal disk centered at $t=1/2$ that contains the  domain of integration of \eqref{integraluno}); also, it satisfies condition (ii), that is, $1/z\notin D$, for $S=\lbrace z\in\CC$, $\Re z<1\rbrace$ (see Fig. 3).

Then, for $\Re z<1$, we can introduce the expansion \eqref{tayloruno} in \eqref{integraluno} and interchange summation and integration to obtain
\begin{equation}\label{expanuno}
{}_2F_1(a,b,c;z)=\left(1-{z\over 2}\right)^{-a}\sum_{n=0}^\infty{(a)_n\over n!}\left({z\over z-2}\right)^n\Phi_n(b,c),
\end{equation}
with 
$$
\Phi_n(b,c):={\Gamma(c)\over\Gamma(b)\Gamma(c-b)}\int_0^1t^{b-1}(1-t)^{c-b-1}(1-2t)^ndt={}_2F_1(-n,b,c;2).
$$

\vskip -3cm
\centerline{\includegraphics[width=4.75in]{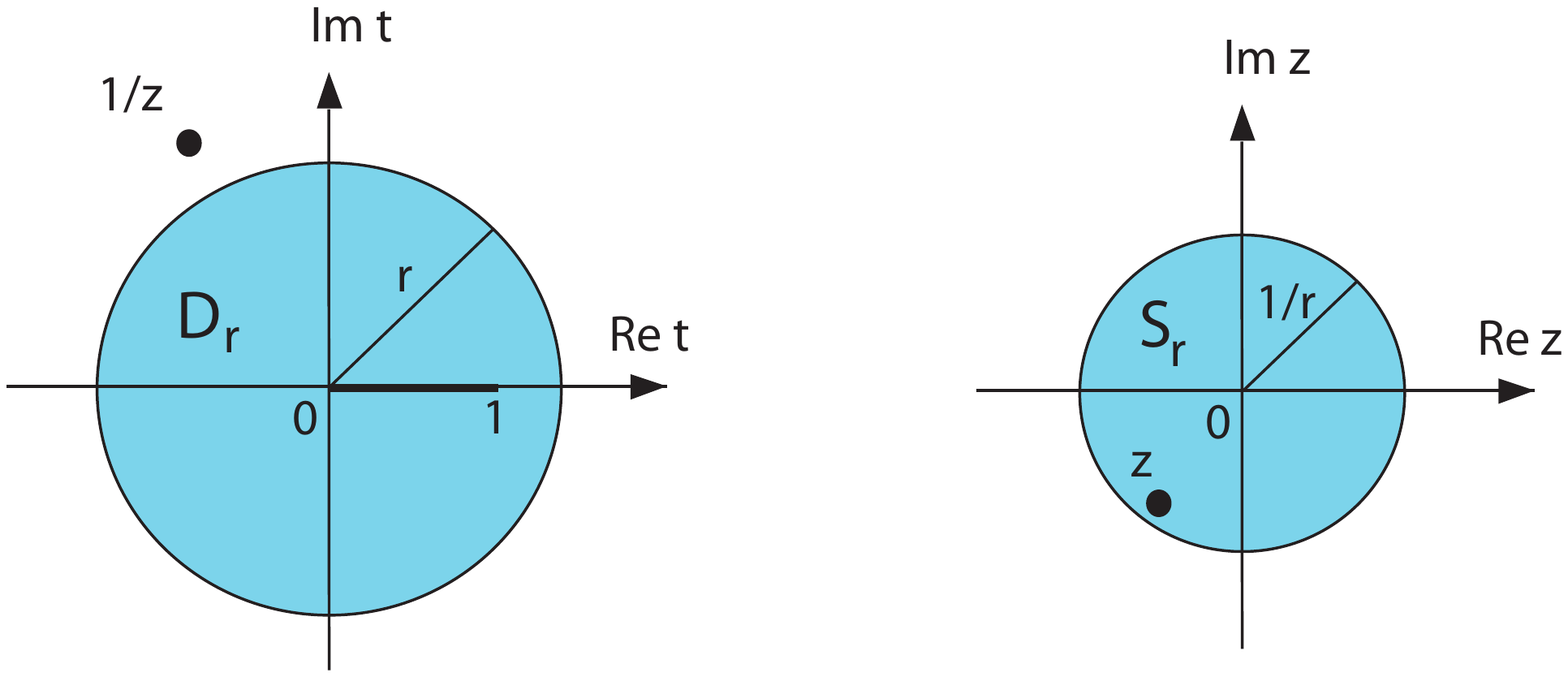}}
\vskip -7.5cm
\centerline{(a)\hskip 7cm (b)}

\noindent
{\bf Figure 2.} The disk of convergence $D_r$ of the Taylor expansion of $f(t)$ at $t=0$ is shown in figure (a) for a certain $r$ ($>1$), and the region $S_r$, inverse of the exterior of $D_r$ is shown in figure (b). The smaller $D_r$ is, the larger $S_r$ is. The smallest possible value of $r$ for which the integration interval $(0,1)\subset D_r$ is $r=1$.
\medskip

Therefore,
\begin{equation}\label{siete}
{}_2F_1(a,b,c;z)=\left(1-{z\over 2}\right)^{-a}\sum_{n=0}^\infty{(a)_n\over n!}\left({z\over z-2}\right)^n{}_2F_1(-n,b,c;2).
\end{equation}
We have $\Phi_0(b,c)=1$, $\Phi_1(b,c)=1-2b/c$ and, for $n=1,2,3,...$, the remaining $\Phi_n(b,c)$ may be obtained from the three-terms recurrence relation \cite[eq. 15.5.11]{nistgauss}
$$
(c+n){}_2F_1(-n-1,b,c;2)+(2b-c){}_2F_1(-n,b,c;2)-n{}_2F_1(1-n,b,c;2)=0.
$$
It is straightforward to show that $\Phi_n(b,c)$ also satisfies the contiguous relation
$$
\Phi_n(b,c)=\Phi_{n-1}(b,c)-{2b\over c}\Phi_{n-1}(b+1,c+1).
$$
The functions $\Phi_n(b,c)$ are polynomials of $b$ and rational functions of $c$. 

%
\centerline{\includegraphics[width=5.5in]{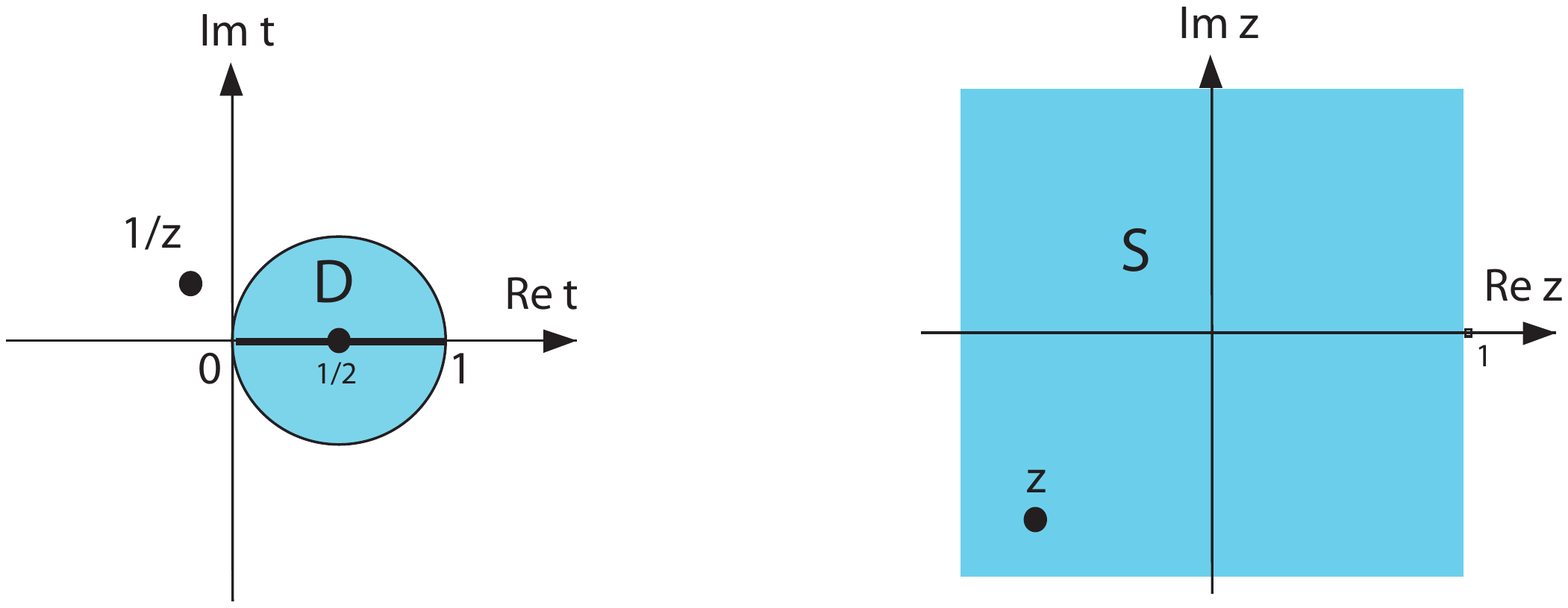}}
\vskip -8.5cm
\centerline{(a)\hskip 8cm (b)}

\noindent
{\bf Figure 3.} The minimal domain of convergence $D$ of the standard Taylor expansion of $f(t)$ at $t=1/2$ containing the interval $(0,1)$ is a disk of radius $1/2$ (figure (a)). The region $S$, inverse of the exterior of $D$ is the region shown in figure (b): $S=\lbrace z\in\CC$, $\Re z<1\rbrace$.
\medskip

The following table shows some numerical experiments comparing the accuracy of B\"u\-hring's expansion and expansion \eqref{siete}. For values of $z$ near the exceptional points $e^{\pm i\pi/3}$, expansion \eqref{siete} is more competitive than B\"uhring's expansion. Away from these points, B\"uhring's expansion becomes more competitive.

\renewcommand{\arraystretch}{1.50}

\bigskip
\centerline{Parameter values: $a=1.2,\ b=2.1, \ c=3, \ z=e^{i\pi/3}, \ z_0=1/2$.} 
\medskip
\begin{tabular}{|l|c|c|c|c|c|}
\hline
\quad\quad\quad\quad$n$ &{0}& {5}& {10}&{15}&{20}\\
\hline
B\"uhring's formula\hfill &{0.263E+2}& {0.879E+1}& {0.103E+1}&{0.955E-1}&{0.803E-2}\\
\hline
Formula \eqref{siete}  &{0.290E+0}&{0.995E-2}&{0.431E-3}&{0.223E-4}&{0.118E-5}\\
\hline
\end{tabular}
\medskip
\centerline{Parameter values: $a=1.2,\ b=2.5, \ c=3, \ z=e^{i\pi/3}, \ z_0=1/2$.} 
\medskip\begin{tabular}{|l|c|c|c|c|c|}
\hline
\quad\quad\quad\quad$n$ &{0}& {5}& {10}&{15}&{20}\\
\hline
B\"uhring's formula\hfill &{0.596E+1}& {0.193E+1}& {0.228E+0}&{0.211E-1}&{0.178E-2}\\
\hline
Formula \eqref{siete}  &{0.467E+0}&{0.228E-1}&{0.126E-3}&{0.734E-4}&{0.437E-5}\\
\hline
\end{tabular}
\medskip
\centerline{Parameter values: $a=1.2,\ b=2.1, \ c=3, \ z=-1, \ z_0=1/2$.} 
\medskip
\begin{tabular}{|l|c|c|c|c|c|}
\hline
\quad\quad\quad\quad$n$ &{0}& {5}& {10}&{15}&{20}\\
\hline
B\"uhring's formula\hfill &{0.155E+2}& {0.330E+0}& {0.248E-2}&{0.148E-4}&{0.796E-7}\\
\hline
Formula \eqref{siete}  &{0.130E+0}&{0.338E-3}&{0.876E-6}&{0.304E-8}&{0.100E-10}\\
\hline
\end{tabular}
\medskip
\centerline{Parameter values: $a=1.2,\ b=2.1, \ c=3, \ z=-1+I, \ z_0=1/2$.} 
\medskip
\begin{tabular}{|l|c|c|c|c|c|}
\hline
\quad\quad\quad\quad$n$ &{0}& {5}& {10}&{15}&{20}\\
\hline
B\"uhring's formula\hfill &{0.114E+2}& {0.972E-1}& {0.291E-3}&{0.690E-6}&{0.149E-8}\\
\hline
Formula \eqref{siete}  &{0.170E+0}&{0.192E-2}&{0.216E-4}&{0.326E-6}&{0.466E-8}\\
\hline
\end{tabular}
\medskip
\centerline{Parameter values: $a=1.2,\ b=2.1, \ c=3.5, \ z=-5, \ z_0=1/2$.} 
\medskip
\begin{tabular}{|l|c|c|c|c|c|}
\hline
\quad\quad\quad\quad$n$ &{0}& {5}& {10}&{15}&{20}\\
\hline
B\"uhring's formula\hfill &{0.475E+1}& {0.407E-3}& {0.619E-8}&{0.852E-13}&{0.1156E-13}\\
\hline
Formula \eqref{siete}  &{0.974E-2}&{0.153E-2}&{0.419E-3}&{0.167E-4}&{0.934E-5}\\
\hline
\end{tabular}
\medskip
\noindent
{\bf Table 1.} The first row represents the number $n$ of terms used in either B\"uhring's expansion or expansion \eqref{siete}. The second row represents the relative error obtained with B\"uhring's approximation. The third row represents the relative error resulting from approximation~\eqref{siete}.
\renewcommand{\arraystretch}{1.0}

\section{An expansion for $2\Re(wz)<1$ with arbitrary $w\in\CC$}

We can generalize the expansion introduced in the above section considering a standard Taylor expansion of the function $f(t)=(1-zt)^{-a}$, not at $t=1/2$, but at a generic point $w=u+iv\in\CC$, $u,v\in\RR$:
\begin{equation}\label{taylorunogen}
f(t)=\sum_{n=0}^\infty{(a)_nz^n\over n!}\left(1-wz\right)^{-a-n}(t-w)^n.
\end{equation}
This expansion satisfies condition (i) for $D=\lbrace t\in\CC$, $\vert t-w\vert<$max$\lbrace\vert w\vert,\vert 1-w\vert\rbrace\rbrace$. It also satisfies condition (ii), that is, $1/z\notin D$, for $S=\lbrace z\in\CC$, $\vert 1-wz\vert>\vert z\vert$max$\lbrace\vert w\vert,\vert 1-w\vert\rbrace\rbrace$. For $u=\Re w\ge 1/2$ the domain $S$ is the semi-plane $S=\lbrace z=x+iy; x,y\in\RR$, $2\Re(wz)=2ux-2vy<1\rbrace$. For $u<1/2$ it is the disk $S=\lbrace z\in\CC$, $\vert z+(1-2u)^{-1}w^*\vert<(1-2u)^{-1}\vert w-1\vert\rbrace$ (see Figure 4).

Then, for $\Re z\in S$, we can introduce the expansion \eqref{taylorunogen} in \eqref{integraluno} and interchange summation and integration to obtain
\begin{equation}\label{expanunogen}
{}_2F_1(a,b,c;z)=\left(1-wz\right)^{-a}\sum_{n=0}^\infty{(a)_n\over n!}\left({wz\over wz-1}\right)^{n}\Phi_n(b,c,w),
\end{equation}
with 
$$
\Phi_n(b,c,w):={\Gamma(c)\over\Gamma(b)\Gamma(c-b)}\int_0^1t^{b-1}(1-t)^{c-b-1}\left(1-{t\over w}\right)^ndt={}_2F_1\left(-n,b,c;{1\over w}\right).
$$

Therefore we have
\begin{equation}\label{diez}
{}_2F_1(a,b,c;z)=\left(1-wz\right)^{-a}\sum_{n=0}^\infty{(a)_n\over n!}\left({wz\over wz-1}\right)^{n}{}_2F_1\left(-n,b,c;{1\over w}\right).
\end{equation}
We have 
$$
\Phi_0(b,c,w)=1,\quad \Phi_1(b,c,w)=1-b/(cw),
$$ 
and, for $n=1,2,3,...$, the remaining $\Phi_n(b,c,w)$ may be obtained from the three-terms recurrence relation \cite[eq. 15.5.11]{nistgauss}
$$
\begin{array}{ll}
\dsp{
(c+n){}_2F_1\left(-n-1,b,c;{1\over w}\right)+\left({b+n\over w}-2n-c\right){}_2F_1\left(-n,b,c;{1\over w}\right)+}\\[8pt]
\dsp{
\quad\quad\quad
n\left(1-{1\over w}\right){}_2F_1\left(1-n,b,c;{1\over w}\right)=0.}
\end{array}
$$
It is straightforward to show that $\Phi_n(b,c)$ also satisfies the contiguous relation
$$
\Phi_n(b,c,w)=\Phi_{n-1}(b,c,w)-{b\over cw}\Phi_{n-1}(b+1,c+1,w).
$$
The functions $\Phi_n(b,c,w)$ are polynomials of $b$ and rational functions of $c$.

\vfill\eject
${}$
\vskip -4cm
\centerline{\includegraphics[width=5.5in]{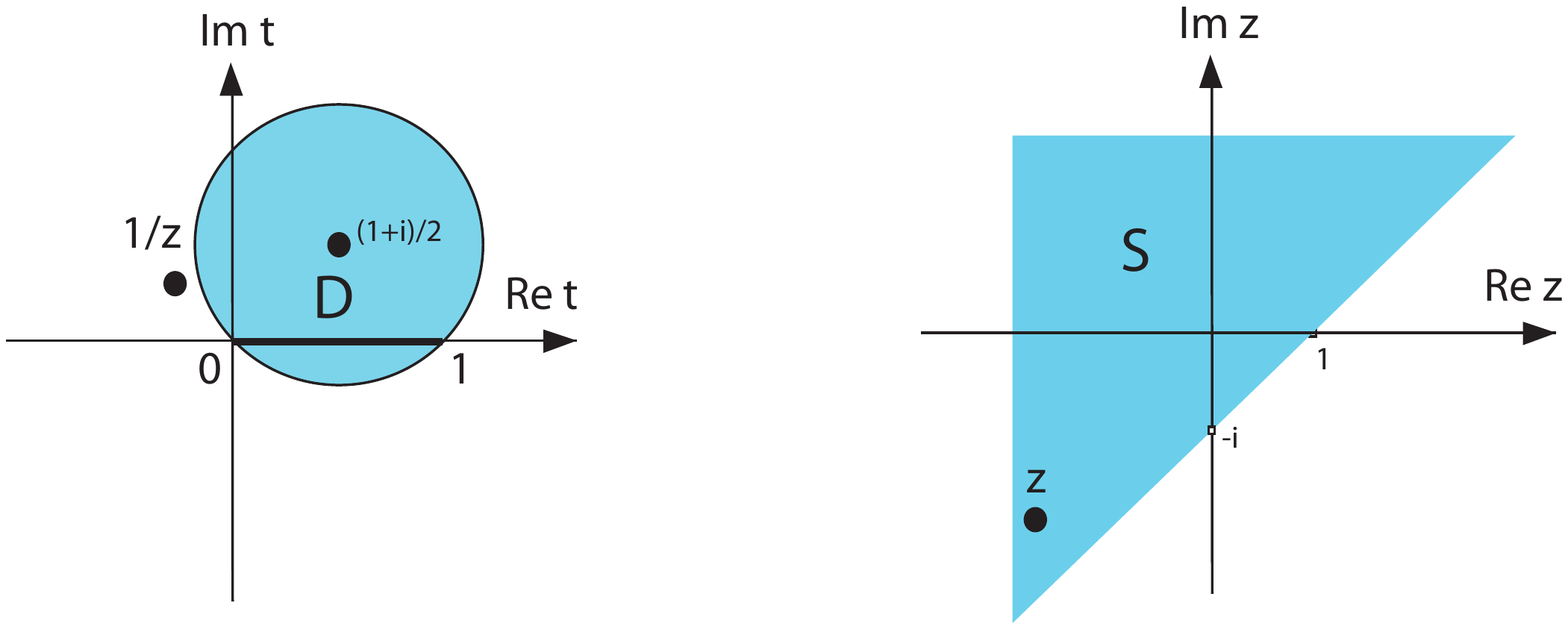}}
\vskip -10.0cm
\centerline{(a) $w=\displaystyle{1+i\over 2}$. \hskip 8cm (b) $w=\displaystyle{1+i\over 2}$.}
\vskip -1.5cm
\medskip
\centerline{\includegraphics[width=4.5in]{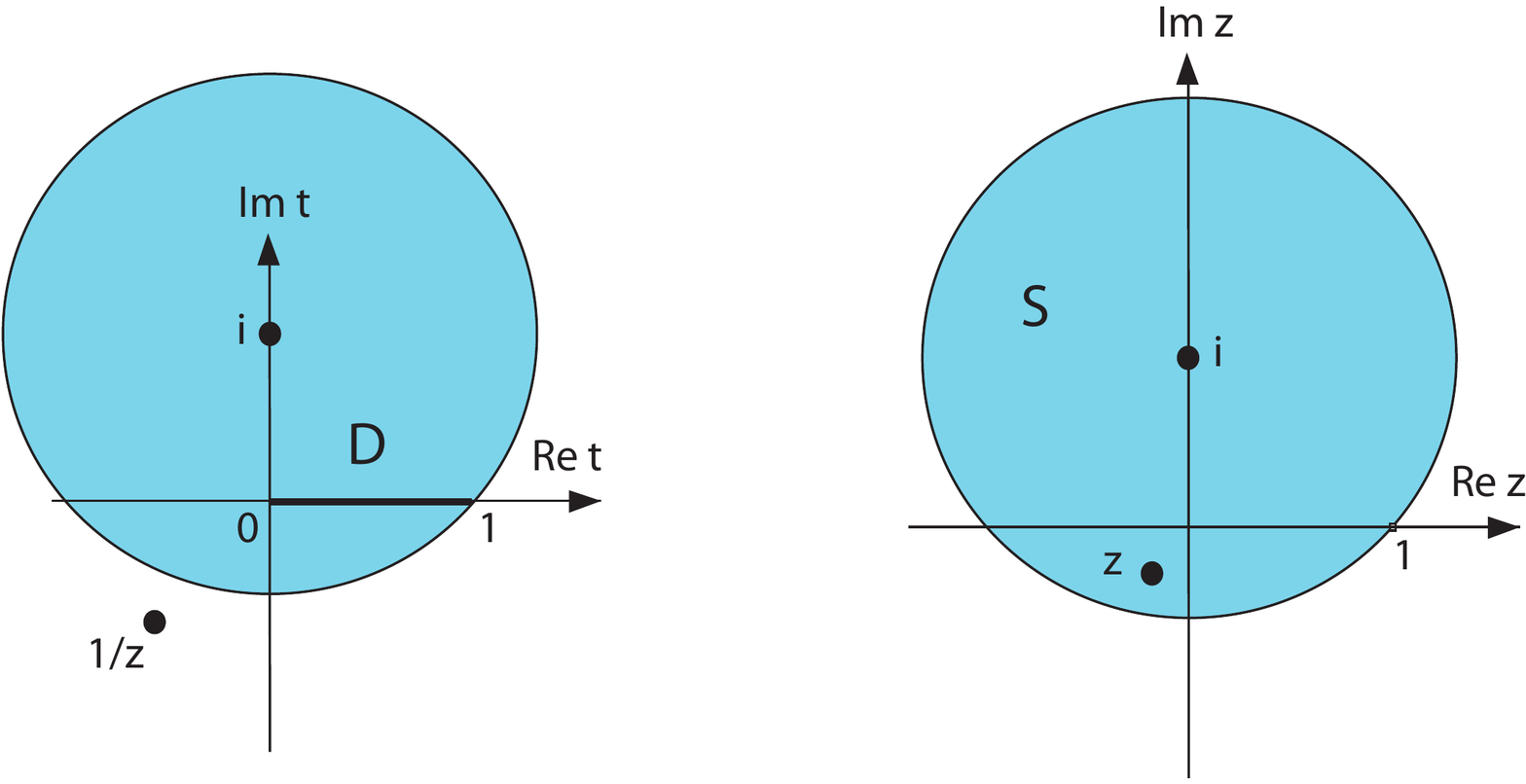}}
\vskip -8cm
\centerline{(c) $w=i$. \hskip 7cm (d) $w=i$.}
\bigskip
\noindent
{\bf Figure 4.} The minimal domain of convergence $D$ of the standard Taylor expansion of $f(t)$ at $t=w$ containing the interval $(0,1)$ is a disk of center at $t=w$ and radius max$\lbrace\vert w\vert,\vert 1-w\vert\rbrace$ (figures (a) and (c)). The region $S$, inverse of the exterior of $D$ is: the half-plane $S=\lbrace z=x+iy; x,y\in\RR$, $1-2\Re(wz)>0\rbrace$ if $\Re w\ge 1/2$ (figure (b)) or the disk of center $w^*/(2\Re w-1)$ and radius $\vert w-1\vert/(1-2\Re w)$ if $\Re w<1/2$ (figure (d)).
\medskip

The following table shows some numerical experiments comparing the accuracy of B\"u\-hring's expansion and expansion \eqref{diez} for $w=(1+i)/2$. For values of $z$ near the exceptional points $e^{\pm i\pi/3}$, expansion \eqref{diez} is more competitive than B\"uhring's expansion. Away from these points, B\"uhring's expansion becomes more competitive.

\renewcommand{\arraystretch}{1.50}

\bigskip
\centerline{Parameter values: $a=1.2,\ b=2.1, \ c=3, \ z=e^{i\pi/3}, \ z_0=1/2$.} 
\medskip
\begin{tabular}{|l|c|c|c|c|c|}
\hline
\quad\quad\quad\quad$n$ &{0}& {5}& {10}&{15}&{20}\\
\hline
B\"uhring's formula\hfill &{0.263E+2}& {0.879E+1}& {0.954E-1}&{0.101E+0}&{0.803E-2}\\
\hline
Formula \eqref{diez}  &{0.408E+0}&{0.606E-2}&{0.156E-3}&{0.476E-5}&{0.150E-6}\\
\hline
\end{tabular}
\medskip
\centerline{Parameter values: $a=1.2,\ b=2.5, \ c=3, \ z=e^{i\pi/3}, \ z_0=1/2$.} 
\medskip
\begin{tabular}{|l|c|c|c|c|c|}
\hline
\quad\quad\quad\quad$n$ &{0}& {5}& {10}&{15}&{20}\\
\hline
B\"uhring's formula\hfill &{0.596E+1}& {0.192E+1}& {0.228E+0}&{0.211E-1}&{0.178E-2}\\
\hline
Formula \eqref{diez}  &{0.480E+0}&{0.127E-1}&{0.408E-3}&{0.138E-4}&{0.477E-6}\\
\hline
\end{tabular}
\medskip
\centerline{Parameter values: $a=1.2,\ b=2.1, \ c=3, \ z=-1, \ z_0=1/2$.} 
\medskip
\begin{tabular}{|l|c|c|c|c|c|}
\hline
\quad\quad\quad\quad$n$ &{0}& {5}& {10}&{15}&{20}\\
\hline
B\"uhring's formula\hfill &{0.154E+2}& {0.330E+0}& {0.248E-2}&{0.148E-4}&{0.796E-7}\\
\hline
Formula \eqref{diez}  &{0.400E+0}&{0.267E-2}&{0.300E-4}&{0.430E-6}&{0.677E-8}\\
\hline
\end{tabular}
\medskip
\centerline{Parameter values: $a=1.2,\ b=2.1, \ c=3, \ z=-1+I, \ z_0=1/2$.} 
\medskip
\begin{tabular}{|l|c|c|c|c|c|}
\hline
\quad\quad\quad\quad$n$ &{0}& {5}& {10}&{15}&{20}\\
\hline
B\"uhring's formula\hfill &{0.114E+2}& {0.972E-1}& {0.291E-3}&{0.690E-6}&{0.149E-8}\\
\hline
Formula \eqref{diez}  &{0.419E+0}&{0.472E-2}&{0.937E-4}&{0.243E-5}&{0.663E-7}\\
\hline
\end{tabular}
\medskip
\centerline{Parameter values: $a=1.2,\ b=2.1, \ c=3.5, \ z=-5, \ z_0=1/2$.} 
\medskip
\begin{tabular}{|l|c|c|c|c|c|}
\hline
\quad\quad\quad\quad$n$ &{0}& {5}& {10}&{15}&{20}\\
\hline
B\"uhring's formula\hfill &{0.475E+1}& {0.407E-3}& {0.619E-8}&{0.852E-13}&{0.156E-14}\\
\hline
Formula \eqref{diez}  &{0.680E+0}&{0.560E-1}&{0.537E-2}&{0.104E-2}&{0.627E-3}\\
\hline
\end{tabular}

\medskip
\noindent
{\bf Table 2.} The first row represents the number $n$ of terms used in either B\"uhring's expansion or expansion \eqref{diez} with $w=(1+i)/2$. The second row represents the relative error obtained with B\"uhring's approximation. The third row represents the relative error resulting from approximation \eqref{diez} with $w=(1+i)/2$.
\renewcommand{\arraystretch}{1.0}

\section{An expansion for $\vert z\vert^2<4\vert1-z\vert$}

As has been pointed out in \cite{nicoester} (in a different context), the use of a multi-point Taylor expansion \cite{nicolopezuno}, \cite{nicolopezdos} with base points in the interval $(0,1)$ is preferable to using a standard Taylor expansion. With a multi-point Taylor expansion we can avoid the singularity $t=1/z$ of $f(t)$ in its domain of convergence in a better way, and, at the same time, include the whole interval $(0,1)$ in its interior (see Fig. 5(a)). 

Therefore, we consider the two-point Taylor expansion of the function $f(t)=(1-zt)^{-a}$ at $t=0$ and $t=1$ \cite{nicolopezuno}:
\begin{equation}\label{taylor}
f(t)=\sum_{n=0}^\infty[A_n(a,z)+B_n(a,z)t]t^n(t-1)^n.
\end{equation}
An explicit formula for the coefficients $A_n(a,z)$ and $B_n(a,z)$ is given in \cite{nicolopezuno}:
$$
A_0(a,z)=1, \hskip 3cm B_0(a,z)=(1-z)^{-a}-1
$$
and, for $n=1,2,3,...$,
$$
A_n(a,z)={1\over n!}\sum_{k=0}^n{(n+k-1)!\over k!(n-k)!}[(-1)^nn-(-1)^kk(1-z)^{k-a-n}](a)_{n-k}z^{n-k}.
$$
$$
B_n(a,z)={1\over n!}\sum_{k=0}^n{(n+k)!\over k!(n-k)!}[(-1)^k(1-z)^{k-a-n}+(-1)^{n+1}](a)_{n-k}z^{n-k}.
$$

Also, a recurrence relation for $A_n(a,z)$ and $B_n(a,z)$ may be obtained by using  the differential equation satisfied by $f(t)$: $(1-zt)f'=azf$. Introducing expansion \eqref{taylor} and
$$
f'(t)=\sum_{n=0}^\infty\lbrace[(2n+1)B_n(a,z)-(n+1)A_{n+1}(a,z)]+(n+1)(2A_{n+1}(a,z)+B_{n+1}(a,z))t\rbrace t^n(t-1)^n,
$$
in the differential equation $(1-zt)f'=azf$, and equating coefficients of $t^n(t-1)^n$ and $t^{n+1}(t-1)^n$ we obtain:
\begin{equation}\label{recu}
\begin{array}{ll}
A_{n+1}(a,z)=&\dsp{{-z(a+2n)A_n(a,z)+[1+n(2-z)]B_n(a,z)\over n+1},} \\[8pt]
B_{n+1}(a,z)=&\dsp{{z(2-z)(a+2n)A_n(a,z)+[z(a+2)+n(6z-z^2-4)-2]B_n(a,z)\over(n+1)(1-z)}. }
\end{array}
\end{equation}

Expansion \eqref{taylor} converges inside a Cassini oval with foci at $t=0$ and $t=1$ and radius $r>0$ of the form $D_r=\lbrace t\in\CC$, $\vert t(t-1)\vert<r\rbrace$. The interval $(0,1)$ is completely contained in this Cassini oval if its middle point $t_0=1/2$ is contained. This happens for $r\ge t_0^2=1/4$ and then, expansion \eqref{taylor} satisfies condition (i) for $r\ge 1/4$. On the other hand, it satisfies condition (ii) if $1/z\notin D_r$ \cite{nicolopezuno}, that is, for any
$$
r<\left\vert{1\over z}\left({1\over z}-1\right)\right\vert.
$$
The smallest $r$ we can take is $r=1/4$ and then, the largest $S_r$ we can choose is (see Fig. 5(b))
$$
S_r=\lbrace z\in\CC;\hskip 2mm\vert z\vert^2<4\vert1-z\vert\rbrace=\lbrace x+iy;x,y\in\RR,\hskip 2mm y^4+(2x^2-16)y^2+[x^4-16x^2+32x-16]<0\rbrace.
$$

Then, for $z\in S_r$, we can introduce the expansion \eqref{taylor} in \eqref{integraluno} and interchange summation and integration to obtain
\begin{equation}\label{expan}
{}_2F_1(a,b,c;z)=\sum_{n=0}^\infty\left[A_n(a,z)\Phi_n(b,c)+B_n(a,z)\Psi_n(b,c)\right],
\end{equation}
with

\vskip -4.5cm
\centerline{\includegraphics[width=5.0in]{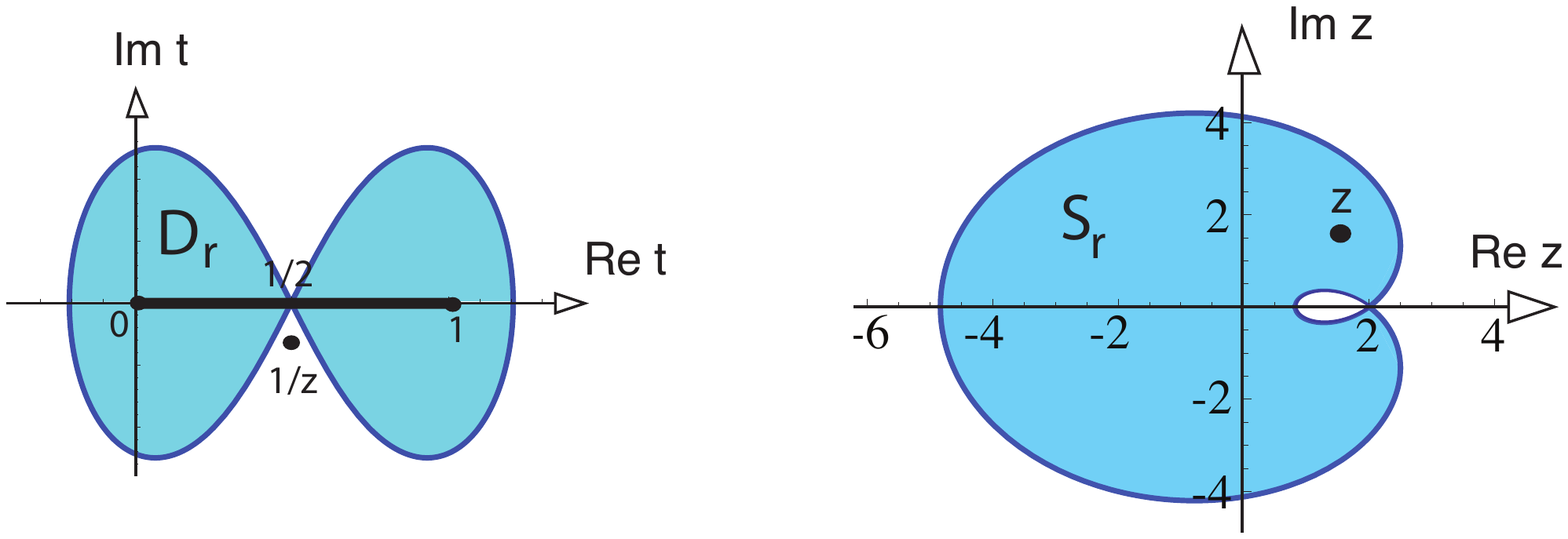}}
\vskip -8cm
\centerline{(a)\hskip 8cm (b)}

\noindent
{\bf Figure 5.}  The minimal domain of convergence $D_r$ of the two-point Taylor expansion of $f(t)$ at $t=0$ and $t=1$ containing the interval $(0,1)$ is a Cassini oval of radius $1/4$ and foci $t=0$ and $t=1$ (figure (a)). The region $S_r$, inverse of the exterior of $D_r$ is the region shown in figure (b): $S_r=\lbrace z\in\CC;\hskip 2mm\vert z\vert^2<4\vert1-z\vert\rbrace$.
\medskip

$$
\Phi_n(b,c):=(-1)^n{\Gamma(c)\over\Gamma(b)\Gamma(c-b)}\int_0^1t^{b+n-1}(1-t)^{n+c-b-1}dt=(-1)^n{(b)_n(c-b)_n\over(c)_{2n}},
$$
$$
\Psi_n(b,c):=(-1)^n{\Gamma(c)\over\Gamma(b)\Gamma(c-b)}\int_0^1t^{b+n}(1-t)^{n+c-b-1}dt=(-1)^n{(b)_{n+1}(c-b)_n\over(c)_{2n+1}}.
$$
Therefore, we have
\begin{equation}\label{expandos}
{}_2F_1(a,b,c;z)=\sum_{n=0}^\infty(-1)^n{(b)_n(c-b)_n\over(c)_{2n+1}}\left[(c+2n)A_n(a,z)+(b+n)B_n(a,z)\right],
\end{equation}
with $A_n(a,z)$ and $B_n(a,z)$  given by the recursion \eqref{recu} and $A_0(a,z)=1$, $B_0(a,z)=(1-z)^{-a}-1$.
This expansion is a series of elementary functions of $z$: a linear combination of $1$ and $(1-z)^{-n-a}$ whose coefficients are polynomials in $z$.

The following table shows some numerical experiments comparing the accuracy of B\"u\-hring's expansion and expansion \eqref{expanunogen}. For values of $z$ near the exceptional points $e^{\pm i\pi/3}$, expansion \eqref{expanunogen} is more competitive than B\"uhring's expansion. Away from these points, B\"uhring's expansion becomes more competitive.

\renewcommand{\arraystretch}{1.50}

\bigskip
\centerline{Parameter values: $a=1.2,\ b=2.1, \ c=3, \ z=-1, \ z_0=1/2$.} 
\medskip
\begin{tabular}{|l|c|c|c|c|c|}
\hline
\quad\quad\quad\quad$n$ &{0}& {5}& {10}&{15}&{20}\\
\hline
B\"uhring's formula\hfill &{0.154E+2}& {0.330E+0}& {0.248E-2}&{0.148E-4}&{0.796E-7}\\
\hline
Formula \eqref{expanunogen}  &{0.112E+0}&{0.242E-5}&{0.630E-10}&{0.143E-14}&{0.408E-15}\\
\hline
\end{tabular}
\medskip
\centerline{Parameter values: $a=1.2,\ b=2.5, \ c=3, \ z=-2, \ z_0=1/2$.} 
\medskip
\begin{tabular}{|l|c|c|c|c|c|}
\hline
\quad\quad\quad\quad$n$ &{0}& {5}& {10}&{15}&{20}\\
\hline
B\"uhring's formula\hfill &{0.181E+1}& {0.294E-2}& {0.175E-5}&{0.805E-9}&{0.339E-12}\\
\hline
Formula \eqref{expanunogen}  &{0.221E+0}&{0.546E-3}&{0.187E-5}&{0.688E-8}&{0.261E-10}\\
\hline
\end{tabular}
\medskip
\centerline{Parameter values: $a=1.2,\ b=2.1, \ c=3, \ z=e^{i\pi/3}, \ z_0=1/2$.} 
\medskip
\begin{tabular}{|l|c|c|c|c|c|}
\hline
\quad\quad\quad\quad$n$ &{0}& {5}& {10}&{15}&{20}\\
\hline
B\"uhring's formula\hfill &{0.263E+2}& {0.879E+1}& {0.955E-1}&{0.101E+0}&{0.803E-2}\\
\hline
Formula \eqref{expanunogen}  &{0.210E+0}&{0.142E-3}&{0.118E-6}&{0.104E-9}&{0.936E-13}\\
\hline
\end{tabular}
\medskip
\centerline{Parameter values: $a=1.2,\ b=2.5, \ c=3, \ z=e^{i\pi/3}, \ z_0=1/2$.} 
\medskip
\begin{tabular}{|l|c|c|c|c|c|}
\hline
\quad\quad\quad\quad$n$ &{0}& {5}& {10}&{15}&{20}\\
\hline
B\"uhring's formula\hfill &{0.596E+1}& {0.193E+1}& {0.228E+0}&{0.211E-1}&{0.178E-2}\\
\hline
Formula \eqref{expanunogen}  &{0.141E+0}&{0.753E-4}&{0.603E-7}&{0.522E-10}&{0.467E-13}\\
\hline
\end{tabular}

\medskip
\noindent
{\bf Table 3.} The first row represents the number $n$ of terms used in either B\"uhring's expansion or expansion \eqref{expanunogen}. The second row represents the relative error obtained with B\"uhring's approximation. The third row represents the relative error resulting from approximation \eqref{expanunogen}.

\renewcommand{\arraystretch}{1.0}

%
\section{An expansion for $\vert z\vert^3<6\sqrt{3}\vert(1-z)(2-z)\vert$}

As has been pointed out in \cite{nicoester}, the use of a three-point Taylor expansion \cite{nicolopezdos} with base points in the interval $(0,1)$ is preferable to the use of a two-point Taylor expansion in order to better avoid the singularity $t=1/z$ of $f(t)$ in its domain of convergence, and, at the same time, to include the whole interval $(0,1)$ in its interior (see Fig. 6(a)). 
Therefore, we consider the three-point Taylor expansion of the function $f(t)=(1-zt)^{-a}$ at $t=0$, $t=1/2$ and $t=1$:
\begin{equation}\label{taylortres}
f(t)=\sum_{n=0}^\infty[A_n(a,z)t+B_n(a,z)t+C_n(a,z)t^2][t(t-1)(t-1/2)]^n,
\end{equation}
with $A_0(a,z)=1$, $B_0(a,z)=4(1-z/2)^{-a}-(1-z)^{-a}-3$ and $C_0(a,z)=2+2(1-z)^{-a}-4(1-z/2)^{-a}$.
We have that \cite{nicoester}\
\begin{equation}\label{taylorprimatres}
f'(t)=\sum_{n=0}^\infty[A'_n(a,z)+B_n'(a,z)t+C'_n(a,z)t^2][t(t-1)(t-1/2)]^n,
\end{equation}
with
$$
\begin{array}{ll}A_n'(a,z):=&{n+1\over 2}A_{n+1}(a,z)+(3n+1)B_n(a,z)+{3n\over 2}C_n(a,z), \\ 
B_n'(a,z):=& (3n+2)C_n(a,z)-3(n+1)A_{n+1}(a,z)-(n+1)B_{n+1}(a,z)-{3(n+1)\over 4}C_{n+1}(a,z), \\
C_n'(a,z):=& 3(n+1)A_{n+1}(a,z)+{3(n+1)\over 2}B_{n+1}(a,z)+{5(n+1)\over 4}C_{n+1}(a,z). 
\end{array}
$$
Introducing \eqref{taylortres} and \eqref{taylorprimatres} into the differential equation $(1-zt)f'=azf$ and equating coefficients of $[t(t-1)(t-1/2)]^n$ we obtain
\begin{equation}\label{recutres}
\begin{array}{ll}
A_{n+1}(a,z)= & {1\over 2(n+1)}\lbrace 2[3n(z-2)-2]B_n(a,z)+4z(3n+a)A_n(a,z)+n(5z-6)C_n(a,z)\rbrace, \\
B_{n+1}(a,z)=&{1\over 2(n+1)(z^2-3z+2)}\lbrace 4z(3n+a)(26z-3z^2-24)A_n(a,z)+\\&
\hskip 0.5cm2[48-4z(18+5a)+6z^2(4+3a)+ \\ & 
\hskip 0.5cm 3n(48-96z+50z^2-3z^3)]B_n(a,z)+[4(20-6z(5+a)+5z^2(2+a))+ \\ & 
\hskip 0.5cm n(264-516z+262z^2-15z^3)]C_n(a,z)\rbrace, \\
C_{n+1}(a,z)=&{1\over (n+1)(z^2-3z+2)}\lbrace 4z(3n+a)(12-12z+z^2)A_n(a,z)+\\&
\hskip 0.5cm 2[2(6(3+a)z-(6+5a)z^2-12)+ \\ & \hskip 0.5cm  3n(z^3-24z^2+48z-24)]B_n(a,z)+[4(2z(9+2a)-3z^2(2+a)-12)+ \\ & 
\hskip 0.5cm n(5z^3-132z^2+276z-144)]C_n(a,z)\rbrace. 
\end{array}
\end{equation}

Expansion \eqref{taylortres} converges inside a Cassini oval with foci at $t=0$, $t=1/2$ and $t=1$ and radius $r>0$ of the form $D_r=\lbrace t\in\CC$, $\vert t(t-1)(t-1/2)\vert<r\rbrace$. The interval $(0,1)$ is completely contained in this Cassini oval if the points $t_0={1\pm\sqrt{3}\over 2}$, at which $R(t):=\vert t(t-1)(t-1/2)\vert$ gets its maximum value, is contained in $D_r$. This happens for $r\ge R(t_0)=(12\sqrt{3})^{-1}$ and then, expansion \eqref{taylortres} satisfies condition (i) for $r\ge (12\sqrt{3})^{-1}$. On the other hand, it satisfies condition (ii) if $1/z\notin D_r$, that is, for any
$$
r<\left\vert{1\over z}\left({1\over z}-1\right)\left({1\over z}-{1\over 2}\right)\right\vert.
$$
The smallest $r$ we can take is $r=(12\sqrt{3})^{-1}$ and then, the largest $S_r$ we can choose is (see Fig. 6(b))
$$
\begin{array}{ll}
S_r
&=\lbrace z\in\CC;\hskip 2mm\vert z\vert^3<6\sqrt{3}\vert(1-z)(2-z)\vert\rbrace\\
&=\lbrace x+iy;x,y\in\RR,\hskip 2mm 108[(1-x)^2+y^2][(2-x)^2+y^2]>(x^2+y^2)^3\rbrace.
\end{array}
$$

Then, for $z\in S_r$, we can introduce the expansion \eqref{taylortres} in \eqref{integraluno} and interchange summation and integration to obtain
\begin{equation}\label{expantres}
\begin{array}{lll}
{}_2F_1(a,b,c;z)&=&\dsp{\sum_{n=0}^\infty(-1)^n\left[A_n(a,z)\Phi_n(b,c)+{b\over c}B_n(a,z)\Phi_n(b+1,c+1)\right.}\\[8pt]
&+&\dsp{\left.{b(b+1)\over c(c+1)}C_n(a,z)\Phi_n(b+2,c+2)\right],}
\end{array}
\end{equation}
with
$$
\begin{array}{ll}
\Phi_n(b,c)&\dsp{:
={\Gamma(c)\over\Gamma(b)\Gamma(c-b)}\int_0^1 t^{n+b-1}(1-t)^{n+c-b-1}(t-1/2)^ndt}\\[8pt]
&\dsp{=(-1)^n{(b)_n(c-b)_n\over 2^n(c)_{2n}} {}_2F_1(-n,b+n;c+2n;2).}
\end{array}
$$
Therefore we have
\begin{equation}\label{ultima}
\begin{array}{ll}
\dsp{{}_2F_1(a,b,c;z)=\sum_{n=0}^\infty\left[A_n(a,z)+{b\over c}B_n(a,z)+{b(b+1)\over c(c+1)}C_n(a,z)\right]\ \times}
\\[8pt]
\quad\quad\dsp{{(b)_n(c-b)_n\over 2^n(c)_{2n}} {}_2F_1(-n,b+n;c+2n;2),}
\end{array}
\end{equation}
with $A_n(a,z)$, $B_n(a,z)$ and $C_n(a,z)$ given by the recursion \eqref{recutres} and $A_0(a,z)=1$, $B_0(a,z)=4(1-z/2)^{-a}-(1-z)^{-a}-3$ and $C_0(a,z)=2+2(1-z)^{-a}-4(1-z/2)^{-a}$.

A recursion relation for $\Phi_n(b,c)$ with respect to $n$ can be obtained by using 
Zeilberger's algorithm\footnote{We thank Dr. Raimundas Vidunas for his  help.} (e.g. via its Maple realization) for hypergeometric functions, 
and it follows that the functions $\Phi_n(b,c)$ satisfy the three-terms recurrence relation
$$
X_n\Phi_{n-1}(b,c)+Y_n\Phi_n(b,c)+Z_n\Phi_{n+1}(b,c)=0,
$$
where
$$
\Phi_0(b,c)=1, \quad\quad\quad \Phi_1(b,c)=-{b(b-c)(2b-c)\over 2c(c+1)(c+2)}.
$$

The coefficients are given by
$$
X_n:=n(-c-2n-5nc-6n^2-4bc+4b^2)(-n+b-c+1)(n+b-1),
$$
$$
Y_n:=2(2b-c)(p_0+p_1n+p_2n^2+p_3n^3),
$$
where
$$
\begin{array}{ll}
p_0:=& 16b(-1+b)(b-c+1)(b-c),\\
p_1:= & -4+21c+40b^2-17c^2-32b^2c+32bc^2-40bc,\\
p_2:= & 24bc+24-24b^2+15c^2-57c,\\
p_3:= & 18(c-2), 
\end{array}
$$

\vfill\eject
${}$
\vskip -4.5cm
\centerline{\includegraphics[width=5in]{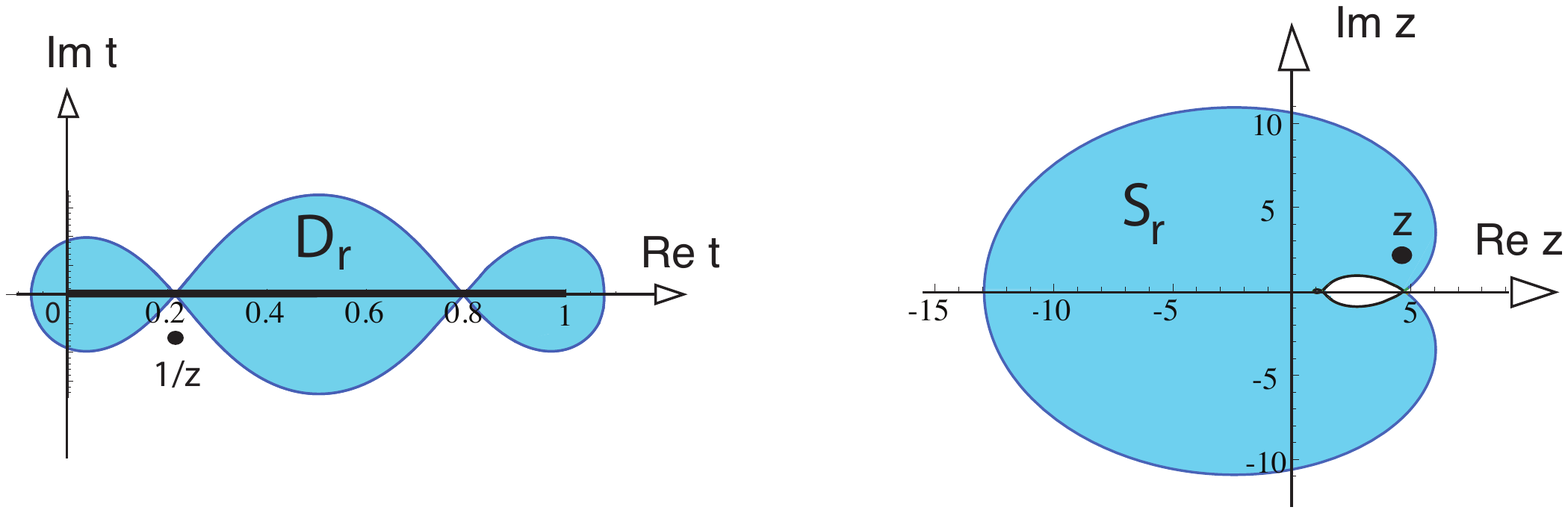}}
\vskip -8.5cm
\centerline{(a)\hskip 8cm (b)}

\noindent
{\bf Figure 6.}  The minimal domain of convergence $D_r$ of the three-points Taylor expansion of $f(t)$ at $t=0$, $t=1/2$ and $t=1$ containing the interval $(0,1)$ is a Cassini oval of radius $(12\sqrt{3})^{-1}$ and foci at $t=0$, $t=1/2$ and $t=1$ (figure (a)). The region $S_r$, inverse of the exterior of $D_r$ is the region shown in figure (b): $S_r=\lbrace z\in\CC;\hskip 2mm\vert z\vert^3<6\sqrt{3}\vert(1-z)(2-z)\vert\rbrace$.
\medskip

and
$$
Z_n:= 16(3n+c))(3n+1+c)(3n+2+c)(-5nc-6n^2+10n+4b^2-4bc+4c-4).
$$

On the other hand, it is straightforward to see that $\Phi_0(b,c)=1$ and that $\Phi_n(b,c)$ satisfies the contiguous relation
$$
\Phi_{n+1}(b,c)={b(b+1)(c-b)\over c(c+1)(c+2)}\Phi_n(b+2,c+3)-{b(c-b)\over 2c(c+1)}\Phi_n(b+1,c+2).
$$

Expansion \eqref{ultima} is a series of elementary functions of $z$: a linear combination of $1$, $(1-z)^{-n-a}$ and $(1-z/2)^{-n-a}$ whose coefficients are polynomials in $z$.

The following table shows some numerical experiments comparing the accuracy of B\"u\-hring's expansion and expansion \eqref{ultima}. For values of $z$ near the exceptional points $e^{\pm i\pi/3}$, expansion \eqref{ultima} is more competitive than B\"uhring's expansion. Away from these points, B\"uhring's expansion becomes more competitive.

\renewcommand{\arraystretch}{1.50}

\bigskip
\centerline{Parameter values: $a=1.2,\ b=2.1, \ c=3, \ z=e^{i\pi/3}, \ z_0=1/2$.} 
\medskip
\begin{tabular}{|l|c|c|c|c|c|}
\hline
\quad\quad\quad\quad$n$ &{0}& {5}& {10}&{15}&{20}\\
\hline
B\"uhring's formula\hfill &{0.263E+2}& {0.177E+2}& {0.879E+1}&{0.253E+1}&{0.103E+1}\\
\hline
Formula \eqref{ultima}  &{0.330E-1}&{0.647E-5}&{0.180E-7}&{0.196E-11}&{0.527E-14}\\
\hline
\end{tabular}
\medskip
\centerline{Parameter values: $a=1.2,\ b=2.5, \ c=3, \ z=e^{i\pi/3}, \ z_0=1/2$.} 
\medskip
\begin{tabular}{|l|c|c|c|c|c|}
\hline
\quad\quad\quad\quad$n$ &{0}& {5}& {10}&{15}&{20}\\
\hline
B\"uhring's formula\hfill &{0.596E+1}& {0.386E+1}& {0.193E+1}&{0.561E+0}&{0.228E-1}\\
\hline
Formula \eqref{ultima}  &{0.351E-1}&{0.496E-5}&{0.137E-7}&{0.184E-11}&{0.523E-14}\\
\hline
\end{tabular}
\medskip
\centerline{Parameter values: $a=1..2,\ b=2.1, \ c=3, \ z=-5, \ z_0=1/2$.} 
\medskip
\begin{tabular}{|l|c|c|c|c|c|}
\hline
\quad\quad\quad\quad$n$ &{0}& {5}& {10}&{15}&{20}\\
\hline
B\"uhring's formula\hfill &{0.841E+0}& {0.165E-2}& {0.206E-4}&{0.236E-7}&{0.238E-9}\\
\hline
Formula \eqref{ultima}  &{0.171E+0}&{0.429E-2}&{0.316E-3}&{0.465E-5}&{0.361E-6}\\
\hline
\end{tabular}
\medskip
\centerline{Parameter values: $a=1.2,\ b=2.01, \ c=3, \ z=-5, \ z_0=1/2$.} 
\medskip
\begin{tabular}{|l|c|c|c|c|c|}
\hline
\quad\quad\quad\quad$n$ &{0}& {5}& {10}&{15}&{20}\\
\hline
B\"uhring's formula\hfill &{0.269E+1}& {0.700E-2}& {0.860E-4}&{0.966E-7}&{0.973E-9}\\
\hline
Formula \eqref{ultima}  &{0.919E-1}&{0.526E-2}&{0.391E-3}&{0.277E-5}&{0.216E-6}\\
\hline
\end{tabular}

\medskip
\noindent
{\bf Table 4.} The first row represents the number $n$ of terms used in either B\"uhring's expansion or expansion \eqref{ultima}. The second row represents the relative error obtained with B\"uhring's approximation. The third row represents the relative error resulting from approximation \eqref{ultima}.
\renewcommand{\arraystretch}{1.0}

\section{Concluding remarks}

In Sections 2 and 3 we have used a standard one-point Taylor expansion of $f(t)=(1-zt)^{-a}$ with the smallest possible convergence region $D$ containing the integration interval $(0,1)$: an expansion at the point $t=1/2$ and convergence radius $r=1/2$ or at any point $t=w$ and convergence radius $r=$max$\lbrace\vert w\vert,\vert w-1\vert\rbrace$. Then, the inverse of the complement of this disk $D$ is the largest possible region that we can obtain with one-point Taylor expansions of $f(t)$: $\Re z<1$ or, in general, the semi-plane $2\Re(zw)<1$ (when $\Re w\ge 1/2$).

In Section 4 we have used a two-point Taylor expansion of $f(t)$, with two base points located in the interval $(0,1)$ and, in Section 5, a three-point Taylor expansion. One may consider the possibility of expanding $f(t)$ at four or more points located in the interval $(0,1)$. In fact, one gets new approximations valid in regions $S_r$ larger than the ones shown in Figs. 5(b) and 6(b); but  the integrals defining the functions $\Phi_n(b,c)$, as well as the recurrences of the coefficients ($A_n(a,z)$, $B_n(a,z)$,...)  become more complicated.

\section{Acknowledgments}
The {\it Direcci\'on General de Ciencia y Tecnolog\'{\i}a} (REF. MTM2010-21037) is acknowledged by its financial support.
NMT acknowledges support from  {\it{Ministerio de Ciencia e Innovaci\'on}}, 
project MTM2009-11686.


\end{document}